\documentclass[11pt]{article}
\usepackage[letterpaper,top=2cm,bottom=2cm,left=3cm,right=3cm,marginparwidth=1.75cm]{geometry}
\usepackage{amsmath}
\usepackage{cite}
\usepackage{url}
\usepackage{graphicx}
\usepackage{color}
\usepackage{siunitx}
\usepackage[utf8]{inputenc}

\usepackage{stmaryrd}
\usepackage{mathtools}
\usepackage{amsfonts}
\usepackage{amsthm}
\usepackage[ruled,linesnumbered]{algorithm2e}

\newcommand{\indices}[2]{\llbracket #1, #2\rrbracket}
\newcommand{\proj}[1]{\mathcal{P}_{#1}}
\newcommand{\A}{\mathbb{A}}
\newcommand{\B}{\mathbb{B}}
\newcommand{\U}{\mathbb{U}}
\newcommand{\V}{\mathbb{V}}
\newcommand{\homog}[1]{\tilde{\mathbb{P}}^{#1}}
\newcommand{\pol}[1]{\mathbb{P}^{#1}}
\newcommand{\newsp}{\mathbb{F}}
\newcommand{\newel}{z}

\newtheorem{theorem}{Theorem}

\title{Generalized Plane Wave quasi-Trefftz spaces for wave propagation in inhomogeneous media %Title (limit the title to a maximum of 17 words)
}

\author%\multauthor
{Ilaria Fontana \hspace{1cm} Lise-Marie Imbert-G\'erard} 
\date{University of Arizona}

\begin{document}

\maketitle
\begin{abstract}
  Partial Differential Equations (PDEs) models for wave propagation in inhomogeneous media are relevant for many applications. We will discuss numerical methods tailored for tackling problems governed by these variable-coefficient PDEs.
  Trefftz methods rely, in broad terms, on the idea of approximating solutions to PDEs via Galerkin methods using basis functions that are exact solutions of the PDE, making explicit use of information about the ambient medium. However, wave propagation in inhomogeneous media is modeled by PDEs with variable coefficients, and in general no exact solutions are available. 
  Quasi-Trefftz methods have been introduced, in the case of the Helmholtz equation, to address this problem: they rely instead on high-order approximate solutions constructed locally. We will discuss basis of Generalized Plane Waves, a particular kind of quasi-Trefftz functions, and how their construction can be related to the construction of polynomial quasi-Trefftz bases.

\end{abstract}
%\keywords{
{\bf keywords: }wave propagation in inhomogeneous media, variable-coefficient partial differential equation, Discontinuous Galerkin method, quasi-Trefftz method, discrete spaces of quasi-Trefftz functions.
%}
%

\section{Introduction}\label{sec:introduction}
From radio-wave propagation in the ionosphere, to underwater acoustics, to heating and probing fusion plasmas, wave propagation in inhomogeneous media has many applications. Generalized Plane Waves (GPWs) were first introduced in \cite{10.1093/imanum/drt030} for a time-harmonic problem of reflectometry in plasma physics, modeled by the so-called O-mode equation, namely:
\begin{equation}
    -\Delta u -\frac{\omega^2}{c^2}\left( 1-\frac{\omega_p^2}{\omega^2} \right) u = 0.
\end{equation}
In this Helmholtz equation with a variable wavenumber, $\omega$ is the frequency; $c$ is the speed of light in vacuum; and $\omega_p$ is the plasma frequency. The latter is proportional to the square of the plasma density, which depends on the position in space. For low density, $\omega_p^2<\omega$, waves can propagate, but at high density, $\omega_p^2>\omega$, the medium is evanescent. The density such that $\omega_p^2=\omega^2$ is called the cut-off density, and for a given plasma density, changing the frequency $\omega$ changes the position of the cut-off in the domain.

Providing a high-order numerical method for simulations of waves bouncing on a cut-off was the original motivation to develop GPW-based quasi-Trefftz methods, taking advantage of a particular type of efficient high-order methods --the so-called Trefftz methods-- while handling variable-coefficient PDEs. A brief overview of these methods is presented in Section \ref{TnQT}. 
Generally speaking, the specificity of these methods is to rely on problem-dependent discrete spaces to represent the unknown: these spaces are designed to retain good approximation properties (see Section \ref{ssec:BAP}) but with as few degrees of freedom as possible.
Previous works \cite{10.1007/s00211-015-0704-y,10.1007/s00211-021-01220-9} have focused on GPW functions for specific problems. The present work is the first effort to present an abstract framework to study GPW spaces, in parallel to the abstract framework for polynomial quasi-Trefftz presented in \cite{QTspaces}. Even though the GPW ansatz was introduced for the scalar wave propagation problem, one may seek a new ansatz of quasi-Trefftz functions for other types of problems, which the abstract framework developed here can handle.

\subsection{Trefftz and quasi-Trefftz methods}
\label{TnQT}

%IMPORTANT notion of "quasi" for us is Taylor-based

%maybe end by introducing the notion of qT property, qT operator and qT space/qT basis

Discontinuous Galerkin (DG) methods are numerical techniques for solving PDE-driven problems, relying on a variational framework and on the definition of discrete spaces composed of discontinuous piecewise polynomials. Trefftz DG methods rely on discrete spaces using functions that are locally tailored to the specific equation under consideration. Specifically, the test and trial functions are exact solutions of the governing PDE on each mesh element. The properties of this high-order approach have been extensively studied particularly in the context of wave problems \cite{Hiptmair:16}. Notably, it has been shown that these techniques achieve a given level of accuracy with a reduced number of degrees of freedom compared to standard polynomial DG approaches. 
However, the main limitation of Trefftz methods is that when the PDE contains variable coefficients or a non-zero right-hand side, local exact solutions are generally unavailable. %, making the Trefftz approach inapplicable.

To address this challenge, quasi-Trefftz methods rely on approximate (rather than exact) local solutions in the sense of a Taylor series truncation. More precisely,
given a PDE defined, for a linear differential operator $\mathcal{L}$, by
\begin{equation}
    \mathcal{L} u = 0,
\end{equation} 
the so-called \emph{quasi-Trefftz property} for a function $\varphi$ is
\begin{equation}\label{eq:QT property}
    [T_q\circ \mathcal{L}] [\varphi] = 0,
\end{equation}
where $T_q$ denotes the Taylor operator up to a suitable degree $q$ and with center $\mathbf{x}_0$, and $T_q\circ \mathcal{L}$ is referred to as the {\it quasi-Trefftz operator}. Discrete spaces of functions with this property are referred to as \emph{quasi-Trefftz spaces}.
%the Taylor expansion up to a suitable degree and with center $\mathbf{x}_T$ of the operator applied to the test and trial functions is zero. 

Various families of quasi-Trefftz functions have been proposed \cite{Imbert-Gerard:24}. A first choice is to construct a subspace of the polynomial space $\pol{p}$. Recent studies have demonstrated the convergence and stability of polynomial quasi-Trefftz methods for some scalar problems, for example, for the diffusion-advection-reaction equation \cite{10.1093/imanum/drae094} or for the Schrodinger equation \cite{10.1007/s10444-024-10108-9}.
% {\color{blue} For polynomial methods like these, as discussed in \cite{Imbert-Gerard:24}, the choice $q=p-\gamma$, where $\gamma$ is the order of the differential operator $\mathcal{L}$, is optimal.}
Different types of GPW quasi-Trefftz functions have been studied in \cite{10.1093/imanum/drt030,10.1007/s00211-015-0704-y,10.1007/s00211-021-01220-9}.
The original idea behind GPWs was to start from classical plane waves, that had been efficiently leveraged in Trefftz methods for problems modeling wave propagation in homogeneous media, and to allow for further flexibility via higher-order terms, either in the phase or in the amplitude. The guiding thread was to retain the oscillating behavior, since it provided high-order approximation properties with less degrees of freedom than classical polynomial spaces. %, while , because .
The main goal of this paper is to investigate the general properties of quasi-Trefftz spaces of phase-based GPWs, i.e., functions $e^P$ where $P$ is a polynomial.

\subsection{Quasi-Trefftz spaces}
% {\color{red}I think we need to make here directly the choice $q=p-2$ once and for all maybe referring to the three-type paper for justification - and clean up this section accordingly}\\
%{\color{violet} \it IF Rem: I tried different possible locations for this statement; I think the best option is either here or at the end of the previous subsection (because this can only be stated in this way after we introduce $q$, $p$, and the fact that GPWs involve a polynomial $P$.} \\
From now on, if $\gamma$ is the order of the differential operator $\mathcal{L}$, we will fix $q=p-\gamma$ in the Taylor operator of the quasi-Trefftz property \eqref{eq:QT property}, since, as discussed in \cite{Imbert-Gerard:24}, this choice is optimal for both polynomials and GPWs.

For polynomial quasi-Trefftz spaces,
the operator $T_{p-\gamma}\circ \mathcal L$
between $\mathbb P$ and $\pol{p-\gamma}$  is linear;
the dimension of the kernel of $T_{p-\gamma}\circ \mathcal L$ is infinite;
the dimension of the kernel of the operator 
\begin{equation}
    \mathcal T_{pol}:=T_{p-\gamma}\circ \mathcal L|_{\pol{p}}
\end{equation}
is finite. In this context, the quasi-Trefftz space is defined as %$\ker\mathcal T_{pol}=\{\varphi \in\pol{p}, \mathcal T_{pol}[\varphi]=0\}$.
\begin{equation}
    \ker\mathcal T_{pol}=\{\varphi \in\pol{p}, \mathcal T_{pol}[\varphi]=0\}.
\end{equation}
A block-triangular structure of the operator $\mathcal T_{pol}$ can be leveraged to construct a quasi-Trefftz basis, and this approach is presented in \cite{QTspaces}.

Generalized Plane Waves are defined starting from a plane wave but have higher-order terms, either in their amplitude or in their phase. Here, we focus on phase-based GPWs. %, i.e., functions of the form $e^P$ where $P$ is a polynomial. 
In order to leverage GPWs, quasi-Trefftz methods rely on the existence of convenient GPW bases.

For phase-based GPW quasi-Trefftz spaces, the situation differs from the polynomial case. 
First, the operator between $\mathbb P$ and $\pol{p-\gamma}$, 
\begin{equation}
    T_{p-\gamma}\circ \mathcal L\circ \exp:P\to [T_q\circ \mathcal L][e^P]
\end{equation}
is not linear %. But 
as in general $\mathcal{L} e^P = (\mathcal{N}(P)) e^P$, where $\mathcal{N}$ is some nonlinear operator. %We then define a quasi-Trefftz operator as
Here, the quasi-Trefftz operator is 
\begin{equation}\label{eq:Tgendef}
    \mathcal T_{gen}:=T_{p-\gamma}\circ \mathcal N|_{\pol{p}}.
\end{equation}
Second, the dimension of 
$\{\varphi= e^P,P\in\mathbb P, \mathcal T_{gen}[P]=0\}$ is infinite, even if $\mathbb P$ is replaced by $\pol{p}$. 
As a consequence, defining the quasi-Trefftz space as %is defined as $\{\varphi= e^P,P\in\mathbb P_p, \mathcal T[P]=0\}$ %, but its dimension 
\begin{equation}\label{eq:GPWQTspace}
   \mathbb Q\mathbb T_p\coloneq \{\varphi= e^P,P\in\pol{p}, \mathcal T_{gen}[P]=0\},
\end{equation}
its dimension is infinite and $\mathcal T_{gen}$ is nonlinear. The goal is then to construct a set of linearly independent GPWs with desired approximation properties, despite the challenge brought by nonlinearity. 

\subsection{Scope}

%A COMPLETER 

This article focuses on a general approach to achieve this goal, providing an abstract framework that can be applied to other choices of ansatz for quasi-Trefftz functions in the future. 

The procedure to construct a GPW quasi-Trefftz function relies on a fundamental assumption. Given a fixed point $\mathbf x_0$, the differential operator of order $\gamma\in\mathbb N$
\begin{equation}
    \mathcal L_c \coloneq\! \sum_{m=0}^\gamma \sum_{ |\mathbf j|=m} c_{\mathbf j}\left( {\mathbf x}  \right) \partial_{{\mathbf x}}^{\mathbf j},
%\end{equation}  
%this assumption states that the operator
%\begin{equation}
\ \text{ and }\ 
\label{eq:defL*}
    \mathcal L_* \coloneq\! \sum_{ |\mathbf j|=\gamma} c_{\mathbf j}\left( {\mathbf x_0}  \right) \partial_{{\mathbf x}}^{\mathbf j}
\end{equation}
this assumption states that $\mathcal L_*$ 
is surjective between spaces of homogeneous polynomials, namely between $\homog{n+\gamma}$ and $\homog{n}$ for any $n\in\mathbb N$.

%SIMILAR TO THE POL CASE
Thanks to the decomposition of polynomial spaces into direct sums of homogeneous polynomial spaces, the quasi-Trefftz operator $\mathcal T_{gen}$ defined in \eqref{eq:Tgendef} can conveniently be decomposed as the sum of:
\begin{itemize}
% \item its restriction to ${\mathbb P}_{\gamma-1}$: $\mathcal T_{gen}|_{{\mathbb P}_{\gamma-1}}$,
% \item $\mathcal T_{gen}-\mathcal T_{gen}|_{{\mathbb P}_{\gamma-1}}$, in turn split into
% \begin{itemize}
%     \item the terms that do not depend on the (polynomial) independent variable,
%     \item $\mathcal L_*$ defined in \eqref{eq:defL*},
%     \item the remainder that will be denoted $\mathcal R$.
% \end{itemize}
    \item the terms that do not depend on the (polynomial) independent variable, $\mathcal{T}_{gen}^{ind}$,
    \item $\mathcal T_{gen}-\mathcal T_{gen}^{ind}$, in turn split into
    \begin{itemize}
        \item $\mathcal L_*$ defined in \eqref{eq:defL*},
        \item the remainder that will be denoted $\mathcal R$.
    \end{itemize}
\end{itemize}
%{\color{violet} \it IF Rem1: I modified this a bit, putting the P-indepedent terms first (ind = independent), and then only $\mathcal{L}_*$ and $\mathcal{R}$. Indeed now there is not need to write an explicit point for $\pol{\gamma-1}$ since it is included in the two abstract theorems as $\newsp$.}\\
%
%{\color{violet} \it IF Rem2: I think we can remove the following paragraph since this will be explained by the second bullet point in the next subsection.}\\

\subsection{Polynomial versus GPWs}

The first description of an abstract framework to describe quasi-Trefftz spaces was first presented in \cite{QTspaces} for the polynomial case.
In the same spirit, this work introduces an abstract framework entailing the essential properties of the quasi-Trefftz operator for the GPW case. To provide perspective, the two cases are compared in the following way:
%The GPW quasi-Trefftz operator being non-linear, 
%The section addresses two important points:
\begin{itemize}
    \item regarding the quasi-Trefftz operator, it is linear in the polynomial case and nonlinear in the GPW case, due to the choice of ansatz; 
    \item regarding the operator splitting, due to the composition of Taylor truncation and a differential operator,
    \begin{itemize}
        \item both in the polynomial and in the GPW operators, the $\mathcal L_*$ terms play a pivotal role,
        \item the remainder $\mathcal{R}$ is linear in the polynomial case and nonlinear in the GPW case,
        \item there is, in the GPW case, a structure corresponding to the triangular structure in the polynomial case;
    \end{itemize} 
    \item regarding the construction of quasi-Trefftz functions, in both cases, the algorithm relies on an iteration over spaces of homogeneous polynomials of increasing degree, with the resolution of a linear system for the operator $\mathcal L_*$ at each step, corresponding only in the polynomial case to a forward substitution procedure;
    \item regarding the number of degrees of freedom in the construction of an individual quasi-Trefftz function, it is the same in both cases, namely $\dim\ker\mathcal L_*|_{\mathbb P^{p}}$;
    %\item regarding the structure of the quasi-Trefftz operator, {\color{red} add some comment about the "lower triangular" structure - which is really coming from the definition of $\mathcal T_{gen}$ but not from any property of $\mathcal L$};
    \item regarding the range of application, in both cases, the sufficient condition on the differential operator $\mathcal L$ is the surjectivity of $\mathcal L_*$ between spaces of homogeneous polynomials.
\end{itemize}

%{\color{red} 
In the polynomial framework, quasi-Trefftz functions are elements of the kernel of a linear operator. In that case, the questions of existence and uniqueness of quasi-Trefftz functions can be answered directly by linear algebra arguments: there always exists quasi-Trefftz functions since the kernel at least contains the polynomial zero; there exists a unique quasi-Trefftz function if the kernel is trivial and the dimension of the space is the dimension of the kernel. 
By contrast, for a nonlinear operator the questions of existence and uniqueness are open and this work provides a general framework to answer them for GPWs: a constructive algorithm proves their existence and the quest for approximation properties leads to prove that the dimension of the space is actually infinite, but convenient finite bases with desired properties can be constructed.
%}

\section{Notation}

In this paper, blackboard bold letters (e.g. $\A$, $\B$) indicate spaces, while callygraphic letters (e.g. $\mathcal{S}$, $\mathcal{R}$) represent operators between spaces. Moreover, we use the notation $\indices{n_1}{n_2}$, $n_1 < n_2$, to denote the set of all integer numbers from $n_1$ to $n_2$. The set of positive integers is denoted $\mathbb N$ and the set of non-negative integers is denoted $\mathbb N_+$. For any $p\in\mathbb N_+$, $\pol{p}$ denotes the space of polynomials of degree at most equal to $p$, and $\homog{p}$ denotes the space of homogeneous polynomials of degree $p$.

\section{Abstract Framework}
\label{sec:AF}

%{\color{violet} \it IF Rem: I changed $\A$ to $\B$ and $x$ to $y$, since now the structure of $\A$ is slightly modified}\\
Given a graded space $\B = \bigoplus_{n=0}^s \B_n$ with $s+1$ layers for some $s\in\mathbb{N}$, each element $y\in \B$ can be uniquely decomposed as the sum of layer elements: 
\begin{equation}\label{eq:decomposition}
    y = \sum_{n=0}^s y_n\ 
    %\quad
    \text{ where }  \forall n\in\indices{0}{s},\  y_n\in \mathbb B_n.
\end{equation}
%where $x_n\in \A_n$, $n\in\indices{0}{s}$.
We now introduce the projection operators
\begin{equation}\label{eq:projection operator}
    \begin{split}
        \forall n\in\indices{0}{s},\ \ 
        \proj{\B_n} \colon \B &\to \B_n, \\
        y &\mapsto y_n.
    \end{split}
    %\qquad \text{for all } n\in\llbracket 0,s\rrbracket.
\end{equation}

\begin{theorem}\label{th:surjectivity}
    Let $\A = \newsp \oplus (\bigoplus_{n=0}^s \A_n)$ and $\B = \bigoplus_{n=0}^s \B_n$ be two graded spaces, %with $s+1$ layers, 
    $s\in\mathbb{N}$, and let $\mathcal{T}\colon \A\to \B$ be a nonlinear operator. Suppose that $\mathcal{T}$ can be decomposed as the sum of two operators $\mathcal{T} = \mathcal{T}_* + \mathcal{R}$ such that
    \begin{enumerate}
        \item 
        \begin{enumerate}
            \item $\mathcal{T}_* \colon \A\to \B$ is linear; \label{it:T linear}
            \item for all $n\in \indices{0}{s}$, $\mathcal{T}_*(\A_n) \subseteq \B_n$; \label{it:T action}
            \item\label{it:surj*} $\mathcal{T}_*$ is surjective; \label{it:T surjective}
            \item $\mathcal{S}_*\colon \B\to \A$ is a (linear) right inverse of $\mathcal{T}_*$; \label{it:S definition}
            \item for all $n\in \indices{0}{s}$, $\mathcal{S}_*(\B_n) \subseteq \A_n$; \label{it:S action}
        \end{enumerate}
        \item 
        \begin{enumerate}
            \item $\mathcal{R}\colon \A\to \B$ is nonlinear; \label{it:R nonlinear}
            \item for all $n\in\indices{0}{s-1}$ %and all $m\in\indices{n}{s}$, 
            $$\mathcal{R}\left(\bigoplus_{k=n}^s \A_k\right) \subseteq \bigoplus_{k=n+1}^s \B_k;$$ \label{it:R action}
            \item $\mathcal{R}(\A_s) = \{0_{\B}\}$; \label{it:R action s}
            \item for all $n\in\indices{1}{s}$ and for all $x \in \A$ with zero $\newsp$-component, 
            $$\proj{\B_n} \left[ \mathcal{R}\left(\sum_{k=0}^s x_k\right) \right] = \proj{\B_n} \left[ \mathcal{R}\left(\sum_{k=0}^{n-1} x_k\right) \right]$$
            % $$\proj{\B_n} \left[ \mathcal{R}\left( {\color{blue} x} \right) \right] = \proj{\B_n} \left[ \mathcal{R}\left( {\color{blue} \newel +} \sum_{k=0}^{n-1} x_k\right) \right]$$
            where %$x = \sum_{k=0}^s x_k$ as in \eqref{eq:decomposition}, and
            $x = 0_{\newsp} + \sum_{k=0}^s x_k$, with %$\newel\in\newsp$ and 
            $x_k\in\A_k$ for all $k\in\indices{0}{s}$,
            and $\proj{\B_n}\colon \B\to \B_n$ is the projection operator defined as in \eqref{eq:projection operator}.
            \label{it:R projections}
        \end{enumerate}
    \end{enumerate}
    Then, $\mathcal{T}$ is surjective, and a right inverse is the operator $\mathcal{S}\colon \B\to \A$ defined by 
    \begin{equation}
        \mathcal{S} = \mathcal{S}_* \circ \sum_{k=0}^s (-\mathcal{R}\circ \mathcal{S}_*)^k
    \end{equation}
\end{theorem}

\begin{proof}
    Given $y\in \B$, and for any $x\in\A$ with zero $\newsp$-component,
    write %To show subjectivity, we need to find $x\in \A$ such that $\mathcal{T}(x) = y$. To show this, we decompose $x$ and $y$ as the sum of layered elements:
    \begin{equation}
        x = \sum_{k=0}^s x_k
        \qquad\text{and}\qquad
        y = \sum_{k=0}^s y_k
    \end{equation}
    where $(x_k,y_k)\in \A_k\times \B_k$ for all $k\in\indices{0}{s}$.
    %Using the decomposition of $\mathcal{T}$, and the properties of $\mathcal{T}_*$ and $\mathcal{R}$ we get
    Then
    \begin{equation}\label{eq:T splitting}
        \begin{split}
            \mathcal{T}(x) =&\ \mathcal{T}_*(x) + \mathcal{R}(x) = \mathcal{T}_*\left(\sum_{k=0}^s x_k\right) + \mathcal{R}\left(\sum_{k=0}^s x_k\right) \\
            =&\ \sum_{k=0}^s \mathcal{T}_*\left(x_k\right) + \sum_{n=0}^s \proj{\B_n} \left[ \mathcal{R}\left(\sum_{k=0}^s x_k\right) \right] \\
            % =&\ {\color{blue} \mathcal{T}_*\left(\newel + \sum_{k=0}^s x_k\right)} + \sum_{n=0}^s \proj{\B_n} \left[ \mathcal{R}\left( {\color{blue} x} \right) \right] \\
            =&\ \sum_{k=0}^s \mathcal{T}_*\left(x_k\right) + \sum_{n=1}^s \proj{\B_n} \left[ \mathcal{R}\left(\sum_{k=0}^{n-1} x_k\right) \right]
        \end{split} 
    \end{equation}
    As a consequence, $\mathcal{T}(x) = y$ if and only if
    \begin{equation}
        \sum_{k=0}^s \mathcal{T}_*\left(x_k\right) + \sum_{k=1}^s \proj{\B_k} \left[ \mathcal{R}\left(\sum_{l=0}^{k-1} x_l\right) \right] = \sum_{k=0}^s y_k
    \end{equation}
    %and this holds if and only if
    if and only if
    \begin{equation}
        % \Leftrightarrow
        \begin{cases}
            \mathcal{T}_* (x_0) = y_0\\
            \mathcal{T}_*(x_k) = y_k - \proj{\B_k} \left[\mathcal{R}\left(\displaystyle\sum_{l=0}^{k-1} x_l\right)\right] &\ \forall k\in\indices{1}{s}
        \end{cases}
    \end{equation}
    For each $k\in\indices{0}{s}$, if the right-hand side is known, a solution $x_k\in \A_k$ exists under hypotheses \ref{it:surj*} and \ref{it:T action}.
    Hence, a suitable solution $x\in\A$ to the system can be constructed by solving for $x_k$ for increasing values of $k$.

    The components of a possible solution $x\in \A$ are:
    \begin{equation}
        \begin{cases}
            \newel = 0_{\newsp} \\
            x_0 = \mathcal{S}_*(y_0)\\
            x_k = \mathcal{S}_*\left( y_k - \proj{\B_k} \left[\mathcal{R}\left(\displaystyle\sum_{l=0}^{k-1} x_l\right)\right] \right) \ \forall k\in\indices{1}{s}
        \end{cases}
    \end{equation}
    Summing all the components, we get
    \begin{equation}
        \begin{split}
            x =&\ %\sum_{k=0}^s x_k = 
            \mathcal{S}_*(y_0) + \sum_{k=1}^s \mathcal{S}_*\left(y_k - \proj{\B_k}\left[\mathcal{R}\left(\sum_{l=0}^{k-1} x_l\right)\right] \right) \\ 
            =&\ \mathcal{S}_* \left(y - \sum_{k=1}^s \proj{\B_k} \left[\mathcal{R}\left(\sum_{l=0}^s x_l\right)\right] \right) \\
            =&\ \mathcal{S}_* \left(y - \sum_{k=0}^s \proj{\B_k} \left[R(x)\right] \right) 
        \end{split}
    \end{equation}
    So $x= \mathcal{S}_*(y-R(x))$.
    Iterating this $l$ times 
    \begin{equation}
        x = \mathcal{S}_* \left(\sum_{k=0}^l (-\mathcal{R}\circ\mathcal{S}_*)^k (y)\right) - \mathcal{S}_*%\circ 
        (-\mathcal{R}\circ\mathcal{S}_*)^l (\mathcal{R}(x))
    \end{equation}
    Moreover, hypothesis \ref{it:S action} together with \ref{it:R action} and \ref{it:R action s} imply that the operator $\mathcal{R}\circ\mathcal{S}_*\colon \B\to \B$ is nilponent, in particular $(\mathcal{R}\circ\mathcal{S}_*)^{s+1} = 0|_{\B\to \B}$. As a consequence, 
    \begin{equation}
        \begin{split}
            x =&\ \mathcal{S}_* \left(\sum_{k=0}^{s+1} (-\mathcal{R}\circ\mathcal{S}_*)^k (y)\right) - \mathcal{S}_*%\circ 
            (-\mathcal{R}\circ\mathcal{S}_*)^{s+1} (\mathcal{R}(x)) \\
            =&\ \mathcal{S}_* \left(\sum_{k=0}^{s} (-\mathcal{R}\circ\mathcal{S}_*)^k (y)\right)
        \end{split}
    \end{equation}
    which concludes the proof.
\end{proof}

The following theorem extends the result of the previous one by presenting an explicit algorithm to construct elements in the counterimage of $y\in\B$, inspired by the forward substitution technique of the polynomial case.

\begin{theorem}\label{th:counterimage}
    Let $\A = \newsp \oplus (\bigoplus_{n=0}^s \A_n)$ and $\B = \bigoplus_{n=0}^s \B_n$ be two graded spaces, %with $s+1$ layers, 
    $s\in\mathbb{N}$, and let $\mathcal{T}\colon \A\to \B$ be a nonlinear operator. Suppose that $\mathcal{T}$ can be decomposed as the sum of two operators $\mathcal{T} = \mathcal{T}_* + \mathcal{R}$ such that properties \ref{it:T linear}--\ref{it:T surjective} and \ref{it:R nonlinear}--\ref{it:R action s} of Theorem \ref{th:surjectivity} hold. Additionally, assume that $\mathcal{T}_*$ and $\mathcal{R}$ satisfies the following conditions:
    \begin{enumerate}
        \item $\mathcal{T}_*(\newsp) = \{0_{\B}\}$\label{it:counterimage 1new}
        \item for all $n\in\indices{0}{s}$, $\A_n = \U_n \oplus \V_n$ and $\mathcal{T}_*|_{\V_n}$ is invertible; \label{it:counterimage 2new}
        \item for all $n\in\indices{0}{s}$, $\mathcal{S}_n\colon\B_n\to \A_n$ is the (linear) right inverse of $\mathcal{T}_*|_{\A_n}$ defined by $\mathcal{S}_n = [\mathcal{T}_*|_{\V_n}]^{-1}$; \label{it:counterimage 3new}
        \item for all $x\in\A$, $\proj{\B_0} \left[ \mathcal{R}\left(x\right) \right] = \proj{\B_0} \left[ \mathcal{R}\left(\newel\right) \right]$, and for all $n\in\indices{1}{s}$,
        \begin{equation*}
            \proj{\B_n} \left[ \mathcal{R}\left(x\right) \right] = \proj{\B_n} \left[ \mathcal{R}\left(\newel + \sum_{k=0}^{n-1} x_k\right) \right]
        \end{equation*}
        where $x = \newel + \sum_{k=0}^s x_k$, with $\newel\in\newsp$ and $x_k\in\A_k$ for all $k\in\indices{0}{s}$. \label{it:counterimage 4new}
    \end{enumerate}
    Then, given $y\in\B$, the following algorithm explicitly computes the components of a generic element $x\in\A$ such that $\mathcal{T}(x) = y$.
    
    \begin{algorithm}[ht]
    \caption{Construction of counterimages}\label{alg:two}
    {\bf Input:} {$\newel\in\newsp$ and $u = \displaystyle\sum_{k=0}^s u_k \in%\U = 
    \displaystyle\bigoplus_{k=0}^s \U_k$} 
    %$u \in\U = \displaystyle\bigoplus_{k=0}^s \U_k$, $u = \displaystyle\sum_{k=0}^s u_k$, where $u_k\in\U_k$ for all $k\in\indices{0}{s}$
    
    $x_0 \gets u_0 + \mathcal{S}_0(y_0 {- \proj{\B_0} 
    \left[\mathcal{R}\left(\newel\right)\right]} - \mathcal{T}_*(u_0))$\;
    \For{$k\in\indices{1}{s}$}{
        $x_k \gets u_k + \mathcal{S}_k \left(y_k - \proj{\B_k} \left[\mathcal{R} \left( {\newel+} \displaystyle\sum_{l=0}^{k-1} x_k\right)\right] - \mathcal{T}_*(u_k) \right)$\;
    }
    {\bf Output: } $x = {\newel +} \displaystyle\sum_{k=0}^s x_k$
    \end{algorithm}
\end{theorem}

\begin{proof}
    We need to check that, if $x\in\A$ is the output of the algorithm, then $\mathcal{T}(x) = y$. 
    %Using the properties of the operators $\mathcal{T}_*$ and $\mathcal{R}$, and in particular the fact that $\mathcal{T}_*\circ\mathcal{S}_k = I_{\B_k}$, we have
    Proceeding as in eq. \eqref{eq:T splitting}, we have
    \begin{equation}\label{eq:T splitting 2}
        \begin{split}
            \mathcal{T}(x) 
            =&\ \mathcal{T}_*(x_0) + \sum_{k=1}^s \mathcal{T}_*(x_k) {+ \proj{\B_0} \left[\mathcal{R}\left(\newel\right)\right]} \\
            &\ + \sum_{k=1}^s \proj{\B_k} \left[ \mathcal{R}\left( {\newel +} \sum_{l=0}^{k-1} x_l\right)\right] \\
        \end{split}
    \end{equation}
    We analyze the first two terms separately. In particular, we use the fact that $\mathcal{T}_*\circ\mathcal{S}_k = I_{\B_k}$ for all $k\in\indices{0}{s}$.
    \begin{equation}\label{eq:T(x0)}
        \begin{split}
            \mathcal{T}_*(x_0) =&\ \mathcal{T}_*(u_0 + \mathcal{S}_0(y_0 {- \proj{\B_0} \left[\mathcal{R}\left(\newel\right)\right]} -\mathcal{T}_*(u_0)) \\
            =&\ \mathcal{T}_*(u_0) + y_0 {- \proj{\B_0} \left[\mathcal{R}\left(\newel\right)\right]} - \mathcal{T}_*(u_0) \\
            =&\ y_0 {- \proj{\B_0} \left[\mathcal{R}\left(\newel\right)\right]}
        \end{split}
    \end{equation}
    and for all $k\in\indices{1}{s}$
    \begin{equation}\label{eq:T(xk)}
        \begin{split}
            \mathcal{T}_*(x_k) %=&\ \mathcal{T}_*\Biggl(u_k + \mathcal{S}_k \Biggl(y_k - \proj{\B_k} \left[\mathcal{R} \left( \sum_{l=0}^{k-1} x_k\right)\right] \\ &\ - \mathcal{T}_*(u_k) \Biggr)\Biggr) \\
            =&\ \mathcal{T}_*(u_k) + y_k - \proj{\B_k} \left[ \mathcal{R}\left(\sum_{l=0}^{k-1} x_l\right)\right] - \mathcal{T}_*(u_k) \\
            =&\ y_k - \proj{\B_k} \left[ \mathcal{R}\left(\sum_{l=0}^{k-1} x_l\right)\right]
        \end{split}
    \end{equation}
    Finally, combining \eqref{eq:T splitting 2}, \eqref{eq:T(x0)}, and \eqref{eq:T(xk)} we conclude that $\mathcal{T}(x) = y$.
\end{proof}

\section{Construction of a GPW function}

%In the polynomial case, quasi-Trefftz functions are elements of the kernel of a linear operator. In that case, the questions of existence and uniqueness of quasi-Trefftz functions can be answered directly by linear algebra arguments: there always exists quasi-Trefftz functions since the kernel at least contains the polynomial zero; there exists a unique quasi-Trefftz function if the kernel is trivial; the 
%By contrast, in the GPW case,

This section aims to illustrate the previous results with practical examples to construct quasi-Trefftz functions. % the previous algorithm can be applied to construct GPW quasi-Trefftz functions, by choosing
Choose $s = p-\gamma$, $\mathbb F=\mathbb P^{\gamma-1}$, $\mathbb A_n=\tilde{\mathbb P}^{n+\gamma}$, $\mathbb B_n=\tilde{\mathbb P}^n$ for all $n\in\indices{0}{p-\gamma}$, so $\mathbb A=\mathbb P^p$ and $\mathbb B=\mathbb P^{p-\gamma}$, together with $\mathcal T_{gen}\coloneq T_{p-\gamma}\circ\mathcal N$, and $\mathcal T_* := \mathcal L_* $ defined in \eqref{eq:defL*}.  

\subsection{Surjectivity} \label{sub:surjectivity} % and initialization} 
%In the quasi-Trefftz literature, the term initialization is related to the choice of degrees of freedom in the construction of an individual quasi-Trefftz function. 

%The sufficient condition for the differential operator $\mathcal L$ is the surjectivity of $\mathcal L_*$ between spaces of homogeneous polynomials

A general discussion on operators $\mathcal L_*$ satisfying the theorem's hypotheses is proposed in\cite{QTspaces}.
%
%Both the Helmholtz and the convected Helmholtz operators are of order $\gamma=2$. 
As a particular case, consider operators $\mathcal L_*$ of order $\gamma=2$ for which $\exists \mathbf j\in(\mathbb N_+)^d$ with $|\mathbf j|=\gamma$ and some component $j_k=\gamma$ such that $c_{\mathbf j}(\mathbf x_0)\neq 0$. They are surjective between spaces of polynomials.
Moreover, regarding the inverse operator $\mathcal S_n$ for any $n\in\mathbb N_0$:
consider %$\mathbb V_n=\{P\in\tilde{\mathbb P}^{n+\gamma},  P = \sum_{|\mathbf j|=n+\gamma, j_k\notin\{0,1\}}\alpha_\mathbf j\mathbf X^{\mathbf j}\}$, 
\begin{equation}
    \mathbb V_n = \left\{P\in\tilde{\mathbb P}^{n+\gamma},  P = \sum_{|\mathbf j|=n+\gamma, j_k\notin\{0,1\}}\alpha_\mathbf j\mathbf X^{\mathbf j}\right\},
\end{equation}
then $\mathcal L_*|_{\mathbb V_n}$ is invertible.
The inverse of this invertible operator satisfies the hypothesis on $\mathcal S_n$ from the previous Theorem \ref{th:counterimage}.
This case includes both the Helmholtz and the convected Helmholtz operators.

\subsection{Helmholtz Equation}

Consider the Helmholtz equation $\mathcal L\varphi = 0$, where
\begin{equation}
 \mathcal L \varphi\coloneq    \Delta \varphi + \kappa^2 \varphi ,
\end{equation}
and $\kappa$ is a variable coefficient. Here, $\mathcal{L} e^P =  \mathcal{N}(P)e^P$ with,
for any $P\in\pol{p}$,
\begin{equation}
  \mathcal{N}(P)\coloneq\Delta P + \lvert \mathbf{\nabla} P\rvert^2 +\kappa^2, \\
\end{equation}
\begin{equation}
    T_{p-2} \circ \mathcal{N} (P) = \Delta P + T_{p-2} \lvert \mathbf{\nabla} P\rvert^2 + T_{p-2}\, \kappa^2, %+ T_{p-2}\, \kappa^2.
\end{equation}
\begin{equation}\label{eq:QT Helmholtz}
     \mathcal T_{gen}  (P) = 0
     \ \Leftrightarrow\ 
    \Delta P + T_{p-2} \lvert \mathbf{\nabla} P\rvert^2 = - T_{p-2}\, \kappa^2,
    % T_{p-2} \circ \mathcal{N} (P) = - T_{p-2}\, \kappa^2
\end{equation}
and we accordingly define the operators %$\mathcal{T}(P) \coloneqq \Delta P + T_{p-2} \lvert \mathbf{\nabla} P\rvert^2$. 
\begin{itemize}
\item[$\bullet$] $    \mathcal{T}(P) \coloneqq \Delta P + T_{p-2} \lvert \mathbf{\nabla} P\rvert^2$,
\item[$\bullet$] $\mathcal T_* (P)  \coloneqq  \Delta P$,$\qquad$ $\bullet$
$\mathcal R (P)  \coloneqq T_{p-2} \lvert \mathbf{\nabla} P\rvert^2$.
\end{itemize}

We now check that $\mathcal{T}_*$ and $\mathcal{R}$ satisfy the desired properties.
Notice that by their definition \ref{it:T linear} and \ref{it:R nonlinear} of Theorem~\ref{th:surjectivity} automatically hold. Moreover,
\begin{itemize}
    \item  {\it [\ref{it:T action} of Th.\ref{th:surjectivity} {and \ref{it:counterimage 1new}. of Th.\ref{th:counterimage}} ]} For all $n\in\indices{0}{p-2}$, $\mathcal{T}_*(\homog{n+2})\subseteq \homog{n}$, {and $\mathcal{T}_*(\pol{1}) = \{0_{\B}\}$.}
    
    \item {\it [\ref{it:T surjective} of Th.\ref{th:surjectivity} and \ref{it:counterimage 2new}.-\ref{it:counterimage 3new}. of Th.~\ref{th:counterimage} (which imply \ref{it:S definition}-\ref{it:S action} of Th.~\ref{th:surjectivity})]} See Subsection \ref{sub:surjectivity}. %$\mathcal{T}_*$ is surjective, because of the surjectivity of the Laplacian operator between spaces of homogeneous polynomials, 
    
    \item {\it [\ref{it:R action} and \ref{it:R action s} of Th.\ref{th:surjectivity}]} Given $n\in\indices{0}{p-2}$ and %$P\in\bigoplus_{k=n}^{p-2} \homog{k+2}$, then %$P=\sum_{k=n}^{p-2} P_{k+2}$ 
    \begin{equation}
        P\in\bigoplus_{k=n}^{p-2} \homog{k+2},\ \text{then}\ 
        P=\sum_{k=n}^{p-2} P_{k+2}
    \end{equation}
    where $P_{k+2}\in\homog{k+2}$ for all $k\in\indices{n}{p-2}$, and
    \begin{equation}\label{eq:action decomposition Helmholtz}
        \begin{split}
            \lvert\mathbf{\nabla} P\rvert^2 =&\ \left\lvert\mathbf{\nabla} \left(\sum_{k=n}^{p-2} P_{k+2}\right)\right\rvert^2 \\
            =&\ \sum_{m=2n+2}^{2p} \sum_{k=-1}^{m-1} \mathbf{\nabla} P_{k+2} \cdot \mathbf{\nabla} P_{m-k}.
        \end{split}
    \end{equation}
    This expression uses the decomposition into the sum of homogeneous polynomials, indeed 
    \begin{equation}
        \sum_{k=-1}^{m-1} \mathbf{\nabla} P_{k+2} \cdot \mathbf{\nabla} P_{m-k} \in\homog{m}
    \end{equation}
    for all $m\in\indices{0}{p-2}$.
    Moreover, the action of $\mathcal{R}$ can be written as
    % \begin{equation}
    %     \begin{split}
    %         \mathcal{R}(P) =&\ T_{p-2} \lvert\mathbf{\nabla} P\rvert^2 \\
    %         =&\ \begin{cases}
    %         \displaystyle\sum_{m=2n+2}^{p-2} \sum_{k=-1}^{m-1} \mathbf{\nabla} P_{k+2} \cdot \mathbf{\nabla} P_{m-k} & \text{if } n\in\indices{-1}{\lfloor p/2\rfloor-2} \\
    %         0 & \text{if } n\in\indices{\lfloor p/2\rfloor-1}{p-2}
    %     \end{cases}
    %     \end{split}
    % \end{equation}
    \begin{equation}
        \mathcal{R}(P) = \sum_{m=2n+2}^{p-2} \sum_{k=-1}^{m-1} \mathbf{\nabla} P_{k+2} \cdot \mathbf{\nabla} P_{m-k}
    \end{equation}
    if $n\in\indices{0}{\lfloor p/2\rfloor-2}$, while %$\mathcal{R}(P) = 0$
    \begin{equation}
        \mathcal{R}(P) = 0
    \end{equation}
    if $n\in\indices{\lfloor p/2\rfloor-1}{p-2}$
    since the lowest degree in \eqref{eq:action decomposition Helmholtz} is $2n+2$.
    As a consequence, for all $n\in\indices{0}{\lfloor p/2\rfloor-2}$, 
    \begin{equation}
        \mathcal{R}\left(\bigoplus_{k=n}^{p-2} \homog{k+2}\right) \subseteq
        \bigoplus_{k=2n+2}^{p-2} \homog{k} \subset \bigoplus_{k=n+1}^{p-2} \homog{k},
    \end{equation}
    and for all $n\in\indices{\lfloor p/2\rfloor-1}{p-2}$, 
    \begin{equation}
        \mathcal{R}\left(\bigoplus_{k=n}^{p-2} \homog{k+2}\right) = \{0_{\B}\}.
    \end{equation}

    \item {\it [{\ref{it:counterimage 4new}. of Th.\ref{th:counterimage} (which implies \ref{it:R projections} of Th.\ref{th:surjectivity})}]} %For all $n\in\indices{0}{p-2}$,
    We have
    \begin{equation}
        \mathcal{R}\left(\sum_{k=-2}^{p-2} P_{k+2}\right) = \sum_{m=0}^{p-2} \sum_{k=-1}^{m-1} \mathbf{\nabla} P_{k+2} \cdot \mathbf{\nabla} P_{m-k}
    \end{equation}
    because $\mathbf{\nabla} P_0 \cdot \mathbf{\nabla} P_{k+2} = 0$ for all $k\in\indices{-2}{p-2}$.
    As a consequence, for all $n\in\indices{0}{p-2}$,
    \begin{equation}
        \begin{split}
            \proj{\homog{n}} \Biggl[\mathcal{R}\Biggl(\sum_{k=-2}^{p-2}& P_{k+2}\Biggr)\Biggr] \\
            =&\ \sum_{k=-1}^{n-1} \mathbf{\nabla} P_{k+2} \cdot \mathbf{\nabla} P_{n-k} \\
            =&\ \proj{\homog{n}} \left[\mathcal{R}\left(\sum_{k=-2}^{n-1} P_{k+2}\right)\right]
        \end{split}
    \end{equation}
\end{itemize}
% {\color{red} maybe update this according to new version of theorems} 
% It is also possible to show that the properties \ref{it:counterimage 2new}. and \ref{it:counterimage 3new}. of Theorem~\ref{th:counterimage} (and consequently \ref{it:S definition} and \ref{it:S action} of Theorem~\ref{th:surjectivity}) hold.  However, for the sake of simplicity, we omit this here, as they are only related to properties of the Laplacian operator.

\subsection{Convected Helmholtz Equation}

We now consider the convected Helmholtz equation for the acoustic potential:
\begin{multline}
    \mathbf{\nabla} \cdot (\rho (\mathbf{\nabla}\varphi - (\mathbf{M}\cdot\mathbf{\nabla}\varphi)\mathbf{M} + {\rm i} \kappa \varphi \mathbf{M}) \\ + \rho (\kappa^2\varphi + {\rm i} \kappa \mathbf{M}\cdot\mathbf{\nabla}\varphi) = 0
\end{multline}
where $\rho$ is the fluid density, $\mathbf{M}$ is the rescaled vector-valued fluid velocity, and $\kappa$ is the (constant) wave number.
Using phase-based GPWs, 
%we get
%\begin{equation}
%    \begin{split}
%        \mathcal{L} e^P %=&\ \mathbf{\nabla} \cdot (\rho (\mathbf{\nabla} e^P - (\mathbf{M}\cdot\mathbf{\nabla} e^P)\mathbf{M} + {\rm i} \kappa  e^P \mathbf{M}) + \rho (\kappa^2 e^P + {\rm i} \kappa \mathbf{M}\cdot\mathbf{\nabla} e^P) \\
%        % =&\ \mathbf{\nabla} \cdot (\rho (\mathbf{\nabla} P - (\mathbf{M}\cdot\mathbf{\nabla} P)\mathbf{M} + {\rm i} \kappa \mathbf{M}) e^P) + \rho (\kappa^2 + {\rm i} \kappa \mathbf{M}\cdot\mathbf{\nabla} P) e^P \\
%        % =&\ [\mathbf{\nabla}\cdot (\rho(\mathbf{\nabla} P - (\mathbf{M}\cdot\mathbf{\nabla} P)\mathbf{M} + {\rm i} \kappa \mathbf{M})) + \rho(\mathbf{\nabla} P - (\mathbf{M}\cdot\mathbf{\nabla} P)\mathbf{M} + {\rm i} \kappa \mathbf{M})\cdot\mathbf{\nabla}P \\
%        % &\ + \rho (\kappa^2 + {\rm i} \kappa \mathbf{M}\cdot\mathbf{\nabla} P)] e^P \\
%        =&\ [\mathbf{\nabla}\rho\cdot \mathbf{\nabla} P + \rho\Delta P - \rho (\mathbf{\nabla}\mathbf{M} \mathbf{\nabla} P) \cdot\mathbf{M} \\ 
 %       &\ - (\mathbf{\nabla} (\rho\mathbf{M}) - 2\rho {\rm i} \kappa) (\mathbf{M}\cdot\mathbf{\nabla} P) \\
%        &\ - \rho(\mathbf{\nabla}(\mathbf{\nabla} P) \mathbf{M}) \cdot\mathbf{M} + \mathbf{\nabla} \cdot(\rho {\rm i}\kappa \mathbf{M}) + \rho \lvert\mathbf{\nabla} P\rvert^2 \\ 
%        &\ - \rho (\mathbf{M}\cdot\mathbf{\nabla} P)^2 + \rho \kappa^2] e^P
%    \end{split} 
%\end{equation}
%Then,
a polynomial $P\in\pol{p}$ satisfies the quasi-Trefftz property if
\begin{equation}
    \begin{split}
        \mathcal{T}(P) \coloneqq&\ T_{p-2}[%\rho\Delta P - \rho(\mathbf{\nabla}(\mathbf{\nabla} P) \mathbf{M}) \cdot\mathbf{M} \\ 
        % &\ + \mathbf{\nabla}\rho\cdot \mathbf{\nabla} P - (\mathbf{\nabla} (\rho\mathbf{M}) - 2\rho {\rm i} \kappa) (\mathbf{M}\cdot\mathbf{\nabla} P) \\ 
        % &\ - \rho (\mathbf{\nabla}\mathbf{M} \mathbf{\nabla} P) \cdot\mathbf{M} + \rho \lvert\mathbf{\nabla} P\rvert^2 - \rho (\mathbf{M}\cdot\mathbf{\nabla} P)^2] \\
        \mathbf{\nabla}\rho\cdot \mathbf{\nabla} P + \rho\Delta P - \rho (\mathbf{\nabla}\mathbf{M} \mathbf{\nabla} P) \cdot\mathbf{M} \\ 
        &\ - (\mathbf{\nabla} (\rho\mathbf{M}) - 2\rho {\rm i} \kappa) (\mathbf{M}\cdot\mathbf{\nabla} P) \\
        &\ - \rho(\mathbf{\nabla}(\mathbf{\nabla} P) \mathbf{M}) \cdot\mathbf{M} + \rho \lvert\mathbf{\nabla} P\rvert^2 \\ 
        &\ - \rho (\mathbf{M}\cdot\mathbf{\nabla} P)^2] \\
        =&\ - T_{p-2} [\mathbf{\nabla} \cdot(\rho {\rm i}\kappa \mathbf{M}) + \rho \kappa^2]
    \end{split}
\end{equation}
Also in this case, $\mathcal{T}\colon\pol{p} \to\pol{p-2}$ is a nonlinear operator that can decomposed into the sum of two operators $\mathcal{T} = \mathcal{T}_* + \mathcal{R}$ defined as follows: if
\begin{equation}
    T_{p-2}\, \rho = \sum_{k=0}^{p-2} \rho_k 
    \qquad\text{and}\qquad
    T_{p-2}\, \mathbf{M} = \sum_{k=0}^{p-2} \mathbf{M}_k
\end{equation}
are the truncated Taylor series %up to order $p-2$ 
of $\rho$ and $\mathbf{M}$, then
\begin{equation}
    \mathcal{T}_* P \coloneqq \rho_0\Delta P - \rho_0(\mathbf{\nabla}(\mathbf{\nabla} P) \mathbf{M}_0) \cdot\mathbf{M}_0,
\end{equation}
and $\mathcal{R}(P) \coloneqq \mathcal{T}(P) - \mathcal{T}_* P$.
% \begin{equation*}
%     \begin{split}
%         \mathcal{R}(P) \coloneqq&\ \mathcal{T}(P) - \mathcal{T}_* P \\
%         =&\ T_{p-2} [\rho\Delta P - \rho_0\Delta P - \rho(\mathbf{\nabla}(\mathbf{\nabla} P) \mathbf{M}) \cdot\mathbf{M} + \rho_0(\mathbf{\nabla}(\mathbf{\nabla} P) \mathbf{M}_0) \cdot\mathbf{M}_0 + \mathbf{\nabla}\rho\cdot \mathbf{\nabla} P \\ 
%         &\ - (\mathbf{\nabla} (\rho\mathbf{M}) - 2\rho {\rm i} \kappa) (\mathbf{M}\cdot\mathbf{\nabla} P) - \rho (\mathbf{\nabla}\mathbf{M} \mathbf{\nabla} P) \cdot\mathbf{M} + \rho \lvert\mathbf{\nabla} P\rvert^2 - \rho (\mathbf{M}\cdot\mathbf{\nabla} P)^2].
%     \end{split}
% \end{equation*}
For this second case, it can be shown using an approach similar to that of the previous section that $\mathcal{T}_*$ and $\mathcal{R}$ satisfy the desired properties. The main difference in this case is the need to include the Taylor expansion components of the variable coefficients $\rho$ and $\mathbf{M}$. For the sake of simplicity, we omit the details.

\section{GPW spaces: basis and approximation properties}

Beyond the existence of GPW quasi-Trefftz functions,  further aspects of GPW spaces are still problem-dependent. %: selection of a finite dimensional space, construction of a basis, approximation properties. 
Various of these aspects have been studied in the literature, see \cite{10.1007/s00211-015-0704-y,Imbert-Gerard:24,10.1007/s00211-021-01220-9}.
In particular, linear independence and approximation properties are closely related in the following way: the former is a by-product of the proof of the latter.

This section revisits pre-existing results in the light of the abstract framework from Section \ref{sec:AF}, highlighting the role played by the surjectivity of the quasi-Trefftz operator in proving approximation properties of quasi-Trefftz spaces.

%Unlike in the polynomial case, where $\mathcal T_{pol}\subset\mathbb P^p$ implies that the dimension of $\ker\mathcal T_{pol}$ is finite, the dimension of $\mathbb Q\mathbb T_p$ defined in \eqref{eq:GPWQTspace} is infinite.

%Here some comments on that for GPW spaces (from work done in previous articles).

\subsection{Candidate for a basis}\label{ssec:cand}

The dimension of the GPW quasi-Trefftz space $\mathbb Q\mathbb T_p$ defined in \eqref{eq:GPWQTspace} is infinite. 
Yet a space of finite dimension, and more particularly its basis, are required in practice for TDG methods. The intuition behind the choice of phase-based GPW ansatz is to seek functions under the form
\begin{equation}
    \varphi(\mathbf x) = \exp \left(\mathrm i \kappa \mathbf d\cdot (\mathbf x-\mathbf x_0) + H.O.T.\right)
\end{equation}
where $H.O.T.$ are higher order polynomial terms while $\mathbf d\in\mathbb C^d$, a unit vector, corresponds to the direction of propagation for classical plane waves.
Hence, in the literature for the Helmholtz equation, families of GPWs were introduced by 
\begin{itemize}
    \item initializing the construction of all functions in a similar fashion by setting
    \begin{itemize}
        \item the linear coefficients as $\mathrm i \kappa \mathbf d$;
        \item  all the other free coefficients to zero;
    \end{itemize}
    \item then for each function picking a distinct $\mathbf d$.
\end{itemize}
For a set of distinct directions $\{\mathbf d_\ell, \mathbf d_\ell\neq\mathbf d_{\tilde{\ell}} \forall \ell\neq\tilde{\ell}\}_{\ell\in I}$ for some index set $I$,
the corresponding GPW family, denoted $\{ \varphi_\ell \}_{\ell\in I}$, is related to a standard plane wave basis, namely
$\{\psi_\ell:\mathbf x\to\exp\mathrm i \kappa \mathbf d_\ell\cdot (\mathbf x-\mathbf x_0)\}_{\ell\in I}$, used in Trefftz methods.
Even though such candidate families are not polynomial,  the polynomial space of their Taylor truncations of degree $p$ are of particular interest for their study. 

\subsection{Plane wave case}

    For reference, given a set of distinct directions $\{\mathbf d_\ell, \mathbf d_\ell\neq\mathbf d_{\tilde{\ell}} \forall \ell\neq\tilde{\ell}\}_{\ell\in I}$ for some index set $I$, consider the polynomial space of Taylor truncations of degree $p$ of the corresponding plane waves, namely $\{ T_p\psi_\ell \}_{\ell\in I}$.
    It is known that its dimension has an upper bound, $\dim\{ T_p\psi_\ell \}_{\ell\in I}\leq D_d$, depending on the dimension $d$: $D_{d=2} = 2p+1$ and $D_{d=3} = (p+1)^2$ 
    \cite{Cessthesis,Moiolathesis} .
    In fact, this is related to trigonometric functions for $d=2$, \cite{10.1007/s00211-021-01220-9},
    and to spherical harmonics for $d=3$, \cite{Imbert-Gerard:24}.
    It is also known that spaces of plane waves with dimension at least equal to $D_d$, that is $\dim I\geq D_d$, have approximation properties of order $p$.

    %Hence plane wave spaces of dimension $D_d$ are used for Trefftz methods.

\subsection{Best approximation properties}
\label{ssec:BAP}
Consider a set of directions $\{\mathbf d_\ell, \mathbf d_\ell\neq\mathbf d_{\tilde{\ell}} \forall \ell\neq\tilde{\ell}\}_{1\leq \ell\leq D_d}$.
For the corresponding candidate family from Section \ref{ssec:cand}, best approximation properties have been studied via the polynomial space of its functions' Taylor truncations of degree $p$, namely $\{ T_p\varphi_\ell \}_{1\leq \ell\leq D_d}$. 
The surjectivity of the quasi-Trefftz operator is fundamental to prove that:
\begin{itemize}
    \item 
    $T_p u \in \{ T_p\varphi_\ell \}_{1\leq \ell\leq D_d}$,
    %\item maybe actually $\{ T_p\varphi_\ell \}_{1\leq \ell\leq D_d} = \{ T_p f, \mathcal T f = 0 \}_{1\leq \ell\leq D_d}$ or something like that
    \item $\dim\{ T_p\varphi_\ell \}_{1\leq \ell\leq D_d} = \dim\{ T_p\psi_\ell \}_{1\leq \ell\leq D_d}$.
\end{itemize}

\begin{theorem}
    Assume $\gamma\in\mathbb N$ and $p\in\mathbb N$ with $p\geq \gamma$. 
If the coefficients of  $\mathcal L$ are in $\mathcal C^{p-\gamma}$, consider a set of directions $\{\mathbf d_\ell, \mathbf d_\ell\neq\mathbf d_{\tilde{\ell}} \forall \ell\neq\tilde{\ell}\}_{1\leq \ell\leq D_d}$ and the corresponding GPW quasi-Trefftz space $\mathbb{QT}_p$. Then,
for any function $u\in\mathcal C^{p+\gamma}$ such that $\mathcal L u=0$ in a neighborhood of $\mathbf x_0$, there exists a quasi-Trefftz function $u_a\in\mathbb{QT}_p$ and  $C\in\mathbb R$ independent of $\mathbf x$, but depending on $u$, $\mathbf x_0$ and $p$ such that, in a neighborhood of $\mathbf x_0$, 
\begin{equation}
|u(\mathbf x)- u_a(\mathbf x)|\leq C |\mathbf x-\mathbf x_0|^{p+1}.
\end{equation}
\end{theorem}

\section{FUTURE WORK}
%maybe not so close: GPW for Maxwell 

This work shines a new light on quasi-Trefftz spaces based on a GPW ansatz. Thanks to an abstract framework, it underlines the distinction between GPW-specific properties and more general quasi-Trefftz properties. 
The latter can be applied to other choices of ansatz for quasi-Trefftz functions in the future, among them the surjectivity of the quasi-Trefftz operator is absolutely fundamental.
The former are mainly related to the ansatz, and therefore to plane waves. Hence for other ansatz we anticipate that they would relate to a relevant reference case.

\section{Acknowledgments}
This material is based upon work supported by the U.S. Department of Energy, Office of Science, Office of Advanced Scientific Computing Research, under Award Number DE-SC0024246.

% For bibtex users:
%\bibliography{fa2025_template}

\bibliographystyle{alpha}
\bibliography{fa2025_template}

% For non bibtex users:
%\begin{thebibliography}{citations}
%\bibitem{Author:00}
%E.~Author.
%\newblock The title of the conference paper.
%\newblock In {\em Proc.\ of the European Society on Vibration
%  }, pages 000--111, Chania, Greece, 2018.
%
%\bibitem{Someone:10}
%A.~Someone, B.~Someone, and C.~Someone.
%\newblock The title of the journal paper.
%\newblock {\em Acta Acust united Ac}, A(B):111--222, 2010.
%
%\bibitem{Someone:04}
%X.~Someone and Y.~Someone.
%\newblock {\em The Title of the Book}.
%\newblock S. Hirzel, Stuttgart, Germany, 2012.
%
%\end{thebibliography}

\end{document}